\documentclass{amsart}
\usepackage{amsmath,amssymb}

\newtheorem{theorem}{Theorem}
\newtheorem{corollary}{Corollary}

\usepackage{epsfig}
\usepackage{url}
\usepackage{amssymb}
\usepackage{amsmath}
\usepackage{amsfonts}

\def\Dbar{\leavevmode\lower.6ex\hbox to 0pt{\hskip-.23ex \accent"16\hss}D}

\def\bj{{\mbox{\bf j}}}

\def\cH{{\mathcal{H}}}

\begin{document}

\title{Generalization of Scarpis's theorem on Hadamard matrices} 

\author {Dragomir {\v{Z}}. {\Dbar}okovi{\'c}}
\address{University of Waterloo, 
Department of Pure Mathematics and Institute for Quantum Computing,
Waterloo, Ontario, N2L 3G1, Canada}
\email{djokovic@uwaterloo.ca}


\date{}

\begin{abstract}
A $\{1,-1\}$-matrix $H$ of order $m$ is a Hadamard matrix 
if $HH^T=mI_m$, where $T$ is the transposition operator and 
$I_m$ the identity matrix of order $m$.
J. Hadamard published his paper \cite{Hadam} on Hadamard matrices 
in 1893. Five years later, Scarpis \cite{Scarpis} showed how 
one can use a Hadamard matrix of order $n=1+p$, 
$p\equiv 3 \pmod{4}$ a prime, to construct a bigger 
Hadamard matrix of order $pn$. In this note we show that 
Scarpis's construction can be extended to the more general 
case where $p$ is replaced by a prime power $q$.
\end{abstract}

\maketitle

\section{Introduction}

We fix some notation which will be used throughout this note.
By $\cH_m$ we denote the set of Hadamard matrices of order $m$.
Let $q\equiv 3 \pmod{4}$ be a prime power and set $n=1+q$.
Let $F_q$ be a finite field of order $q$.

Given a bijection $\alpha:\{1,2,\ldots,q\}\to F_q$, 
we shall construct a map 
\begin{equation*}
\varphi_{q,\alpha}:\cH_n\to\cH_{qn}.
\end{equation*}
Consequently, the following theorem holds.

\begin{theorem} \label{glavna}
Let  $q\equiv 3 \pmod{4}$ be a prime power.
If there exists a Hadamard matrix of order $n=1+q$ then there 
exists also a Hadamard matrix of order $qn$. 
\end{theorem}

In the special case, where $q$ is a prime, this theorem was 
proved by Scarpis \cite{Scarpis}. For a nice and short 
description of the original Scarpis's construction see 
\cite{blog}.

By the well known theorem of Paley \cite{Paley}, the hypothesis of the above theorem 
is always satisfied. Thus we have 

\begin{corollary}
If $q \equiv 3 \pmod{4}$ is a prime power, then there exists a Hadamard matrix of order $q(1+q)$.
\end{corollary}

We shall describe a procedure whose input is a Hadamard matrix $A=[a_{i,j}]$ 
of order $n=1+q$ and output a Hadamard matrix $B=\varphi_{q,\alpha}(A)$ of order $qn$. 
For convenience, we set $\alpha_i=\alpha(i)$. 

\section{Construction of $B$}

Step 1: If $a_{1,1}=-1$ then replace $A$ by $-A$. From now on 
$a_{1,1}=1$.

Step 2: For each $i\in\{2,3,\ldots,n\}$ do the following:
if $a_{i,1}=-1$ then multiply the row $i$ of $A$ 
by $-1$, and if $a_{1,i}=-1$ then multiply the column $i$ of 
$A$ by $-1$. The resulting matrix $A$ is independent of 
the order in which these operations are performed.

Note that $A$ is now normalized, i.e., $a_{i,1}=a_{1,i}=1$ for each $i$. 
Denote by $C$ its \emph{core}, i.e., the submatrix of $A$ 
obtained by deleting the first row and the first column of $A$.
For $i\in\{1,2,\ldots,q\}$, we denote by $c_i$ the row $i$ 
of $C$. For convenience, we also set $c(\alpha_i)=c_i$.

The tensor product $X\otimes Y$ of two matrices 
$X=(x_{i,j})$ and $Y$ is the block matrix $[x_{i,j}Y]$. 

Let $\bj$ be the row vector of length $q$ all of whose entries
are $1$. We view $\bj$ also as a $1\times q$ matrix.

Step 3: We partition $B$ into $n$ blocks of size $q\times qn$:
\begin{equation*}
B=\left[ \begin{array}{c} B_0 \\ B_1 \\  \vdots \\
B_{q} \end{array} \right].
\end{equation*}
We set $B_{0}=A'\otimes\bj$ where $A'$ is the submatrix of 
$A$ obtained by deleting the first row of $A$.

Step 4: For $r\in\{1,2,\ldots,q\}$, we partition $B_r$ into 
$n$ blocks of size $q\times q$:
\begin{equation*}
B_r=\left[ B_{r,0} ~ B_{r,1} ~ \cdots ~ B_{r,q} \right].
\end{equation*}
We set $B_{r,0}=\bj^T \otimes c_r$.

It remains to define the blocks $B_{r,i}$ for 
$\{r,i\}\subseteq\{1,2,\ldots,q\}$.

Step 5: For $\{r,i\}\subseteq\{1,2,\ldots,q\}$, we define 
$B_{r,i}$ by specifying that its row $k$ is 
$c(\alpha_i \alpha_r + \alpha_k)$. Thus  
$B_{r,i}=P_{r,i}C$ where $P_{r,i}$ is a permutation matrix. 

This completes the definition of $B$. 

It remains to prove that $B$ is a Hadamard matrix. \\

\section{Proof that $B$ is a Hadamard matrix}

As $B$ is a square $\{1,-1\}$-matrix of order $qn$, it suffices 
to prove that the dot product of any pair of rows of $B$ is 0. 
There are four cases to consider.

(i) Two distinct rows of $B_0$. They are orthogonal because two 
distinct rows of $A'$ are orthogonal.

(ii) Two distinct rows of $B_r$, $r\ne0$. 
Since $A$ is normalized, the dot product $c_r\cdot c_s$  
is $q$ when $r=s$, and $-1$ otherwise.  Hence, the same is true for each of the blocks $B_{r,i}$ for $i\ne0$.
On the other hand, the dot product of any two rows of $B_{r,0}$ 
is $q$. It follows that the dot product of any pair of rows of $B_r$ is 0.

(iii) A row of $B_0$ and a row of $B_s$, $s\ne0$. 

The row $k$ of $B_0$ is $[~ \bj ~\quad c_k \otimes \bj ~]$ and 
the row $l$ of $B_s$ is 
$$
[~ c(\alpha_s)~~
c(\alpha_s\alpha_l+\alpha_1)~~c(\alpha_s\alpha_l+\alpha_2)~~ 
\cdots ~~ c(\alpha_s\alpha_l+\alpha_q)~].
$$
Since all row sums of $C$ are $-1$, it follows that the dot 
product of the two rows above is 0.

(iv) A row of $B_r$ and a row of $B_s$, $0<r<s$.

The dot product of the row $k$ of $B_r$ and the row $l$ of 
$B_s$ is
$$
c(\alpha_r)\cdot c(\alpha_s)+\sum_{i=1}^q 
c(\alpha_i\alpha_r+\alpha_k)\cdot c(\alpha_i\alpha_s+\alpha_l).
$$
Note that 
$c(\alpha_i\alpha_r+\alpha_k)\cdot c(\alpha_i\alpha_s+\alpha_l)$ 
is equal to $-1$ for all $i$ except that it is equal to $q$ for the unique $i\in\{1,2,\ldots,q\}$ for which
$\alpha_i\alpha_r+\alpha_k=\alpha_i\alpha_s+\alpha_l$. 
Since also $c(\alpha_r)\cdot c(\alpha_s)=-1$, it follows that 
the rows of $B_r$ are orthogonal to the rows of $B_s$.

We have shown that $B\in\cH_{qn}$. 
This completes our construction of $\varphi_{q,\alpha}$. \\

The smallest $q$ which satisfies the condition of 
Theorem \ref{glavna} but is not a prime (so Scarpis's 
theorem does not apply) is $q=27$. It gives a Hadamard 
matrix of order $4\cdot189=756$.

We conclude with an open problem: 
Find an analog of our procedure which uses prime powers $q\equiv 1 \pmod{4}$.

The author acknowledges generous support by NSERC.

\end{document}